\documentclass{amsart}
\usepackage{amssymb}
\usepackage{mathrsfs}
\usepackage{euscript}
\usepackage{lipsum}

\usepackage{enumerate}

\newcommand*{\abs}[1]{\left\lvert#1\right\rvert}   
\newcommand*{\set}[1]{\left\{#1\right\}}           
\newcommand*{\brs}[1]{\left(#1\right)}             
\newcommand*{\norm}[1]{\left\Vert#1\right\Vert}    
\newcommand*{\hilb}{\mathcal H}                     

\newcommand*{\mbC}{{\mathbb C}}
\newcommand*{\mbR}{{\mathbb R}}

\newcommand*{\Tr}{\operatorname{Tr}\,}        

    \newtheorem{thm}{Theorem}                     [section]
    \newtheorem{thm*}{Theorem}
    \newtheorem{prop}[thm]{Proposition}
    \newtheorem{lemma}[thm]{Lemma}
    \newtheorem{cor}[thm]{Corollary}

    \newtheorem{lemma*}{Lemma}    

    \newtheorem{rems*}{Remark}   

\newcommand{\ndef}{\newcommand*}
\def\rndef{\renewcommand}

\ndef{\myaddress}[1]{\begin{center} \it\small #1 \end{center}}

\ndef{\clA}{{\mathcal A}} \ndef{\rmA}{{\mathrm A}} \ndef{\mbA}{{\mathbb A}} \ndef{\bfA}{{\mathbf A}} \ndef{\euA}{{\EuScript A}} \ndef{\frA}{{\mathfrak A}}
\ndef{\clB}{{\mathcal B}} \ndef{\rmB}{{\mathrm B}} \ndef{\mbB}{{\mathbb B}} \ndef{\bfB}{{\mathbf B}} \ndef{\euB}{{\EuScript B}} \ndef{\frB}{{\mathfrak B}}
\ndef{\clC}{{\mathcal C}} \ndef{\rmC}{{\mathrm C}}                          \ndef{\bfC}{{\mathbf C}} \ndef{\euC}{{\EuScript C}} \ndef{\frC}{{\mathfrak C}}
\ndef{\clD}{{\mathcal D}} \ndef{\rmD}{{\mathrm D}} \ndef{\mbD}{{\mathbb D}} \ndef{\bfD}{{\mathbf D}} \ndef{\euD}{{\EuScript D}} \ndef{\frD}{{\mathfrak D}}
\ndef{\clE}{{\mathcal E}} \ndef{\rmE}{{\mathrm E}} \ndef{\mbE}{{\mathbb E}} \ndef{\bfE}{{\mathbf E}} \ndef{\euE}{{\EuScript E}} \ndef{\frE}{{\mathfrak E}}
\ndef{\clF}{{\mathcal F}} \ndef{\rmF}{{\mathrm F}} \ndef{\mbF}{{\mathbb F}} \ndef{\bfF}{{\mathbf F}} \ndef{\euF}{{\EuScript F}} \ndef{\frF}{{\mathfrak F}}
\ndef{\clG}{{\mathcal G}} \ndef{\rmG}{{\mathrm G}} \ndef{\mbG}{{\mathbb G}} \ndef{\bfG}{{\mathbf G}} \ndef{\euG}{{\EuScript G}} \ndef{\frG}{{\mathfrak G}}
\ndef{\clH}{{\mathcal H}} \ndef{\rmH}{{\mathrm H}} \ndef{\mbH}{{\mathbb H}} \ndef{\bfH}{{\mathbf H}} \ndef{\euH}{{\EuScript H}} \ndef{\frH}{{\mathfrak H}}
\ndef{\clI}{{\mathcal I}} \ndef{\rmI}{{\mathrm I}} \ndef{\mbI}{{\mathbb I}} \ndef{\bfI}{{\mathbf I}} \ndef{\euI}{{\EuScript I}} \ndef{\frI}{{\mathfrak I}}
\ndef{\clJ}{{\mathcal J}} \ndef{\rmJ}{{\mathrm J}} \ndef{\mbJ}{{\mathbb J}} \ndef{\bfJ}{{\mathbf J}} \ndef{\euJ}{{\EuScript J}} \ndef{\frJ}{{\mathfrak J}}
\ndef{\clK}{{\mathcal K}} \ndef{\rmK}{{\mathrm K}} \ndef{\mbK}{{\mathbb K}} \ndef{\bfK}{{\mathbf K}} \ndef{\euK}{{\EuScript K}} \ndef{\frK}{{\mathfrak K}}
\ndef{\clL}{{\mathcal L}} \ndef{\rmL}{{\mathrm L}} \ndef{\mbL}{{\mathbb L}} \ndef{\bfL}{{\mathbf L}} \ndef{\euL}{{\EuScript L}} \ndef{\frL}{{\mathfrak L}}
\ndef{\clM}{{\mathcal M}} \ndef{\rmM}{{\mathrm M}} \ndef{\mbM}{{\mathbb M}} \ndef{\bfM}{{\mathbf M}} \ndef{\euM}{{\EuScript M}} \ndef{\frM}{{\mathfrak M}}
\ndef{\clN}{{\mathcal N}} \ndef{\rmN}{{\mathrm N}}                          \ndef{\bfN}{{\mathbf N}} \ndef{\euN}{{\EuScript N}} \ndef{\frN}{{\mathfrak N}}
\ndef{\clO}{{\mathcal O}} \ndef{\rmO}{{\mathrm O}} \ndef{\mbO}{{\mathbb O}} \ndef{\bfO}{{\mathbf O}} \ndef{\euO}{{\EuScript O}} \ndef{\frO}{{\mathfrak O}}
\ndef{\clP}{{\mathcal P}} \ndef{\rmP}{{\mathrm P}} \ndef{\mbP}{{\mathbb P}} \ndef{\bfP}{{\mathbf P}} \ndef{\euP}{{\EuScript P}} \ndef{\frP}{{\mathfrak P}}
\ndef{\clQ}{{\mathcal Q}} \ndef{\rmQ}{{\mathrm Q}}                          \ndef{\bfQ}{{\mathbf Q}} \ndef{\euQ}{{\EuScript Q}} \ndef{\frQ}{{\mathfrak Q}}
\ndef{\clR}{{\mathcal R}} \ndef{\rmR}{{\mathrm R}}                          \ndef{\bfR}{{\mathbf R}} \ndef{\euR}{{\EuScript R}} \ndef{\frR}{{\mathfrak R}}  
\ndef{\clS}{{\mathcal S}} \ndef{\rmS}{{\mathrm S}} \ndef{\mbS}{{\mathbb S}} \ndef{\bfS}{{\mathbf S}} \ndef{\euS}{{\EuScript S}} \ndef{\frS}{{\mathfrak S}}
\ndef{\clT}{{\mathcal T}} \ndef{\rmT}{{\mathrm T}} \ndef{\mbT}{{\mathbb T}} \ndef{\bfT}{{\mathbf T}} \ndef{\euT}{{\EuScript T}} \ndef{\frT}{{\mathfrak T}}
\ndef{\clU}{{\mathcal U}} \ndef{\rmU}{{\mathrm U}} \ndef{\mbU}{{\mathbb U}} \ndef{\bfU}{{\mathbf U}} \ndef{\euU}{{\EuScript U}} \ndef{\frU}{{\mathfrak U}}
\ndef{\clV}{{\mathcal V}} \ndef{\rmV}{{\mathrm V}} \ndef{\mbV}{{\mathbb V}} \ndef{\bfV}{{\mathbf V}} \ndef{\euV}{{\EuScript V}} \ndef{\frV}{{\mathfrak V}}
\ndef{\clW}{{\mathcal W}} \ndef{\rmW}{{\mathrm W}} \ndef{\mbW}{{\mathbb W}} \ndef{\bfW}{{\mathbf W}} \ndef{\euW}{{\EuScript W}} \ndef{\frW}{{\mathfrak W}}
\ndef{\clX}{{\mathcal X}} \ndef{\rmX}{{\mathrm X}} \ndef{\mbX}{{\mathbb X}} \ndef{\bfX}{{\mathbf X}} \ndef{\euX}{{\EuScript X}} \ndef{\frX}{{\mathfrak X}}
\ndef{\clY}{{\mathcal Y}} \ndef{\rmY}{{\mathrm Y}} \ndef{\mbY}{{\mathbb Y}} \ndef{\bfY}{{\mathbf Y}} \ndef{\euY}{{\EuScript Y}} \ndef{\frY}{{\mathfrak Y}}
\ndef{\clZ}{{\mathcal Z}} \ndef{\rmZ}{{\mathrm Z}}                          \ndef{\bfZ}{{\mathbf Z}} \ndef{\euZ}{{\EuScript Z}} \ndef{\frZ}{{\mathfrak Z}}

\ndef{\tA}{{\widetilde A}} \ndef{\tcA}{{\widetilde\clA}} \ndef{\ttcA}{\widetilde{\tcA}} \ndef{\sfA}{{\textsf A}} \ndef{\ttA}{\widetilde{\tA}} \ndef{\dzA}{{A^\sharp}}
\ndef{\tB}{{\widetilde B}} \ndef{\tcB}{{\widetilde\clB}} \ndef{\ttcB}{\widetilde{\tcB}} \ndef{\sfB}{{\textsf B}} \ndef{\ttB}{\widetilde{\tB}} \ndef{\dzB}{{B^\sharp}}
\ndef{\tC}{{\widetilde C}} \ndef{\tcC}{{\widetilde\clC}} \ndef{\ttcC}{\widetilde{\tcC}} \ndef{\sfC}{{\textsf C}} \ndef{\ttC}{\widetilde{\tC}} \ndef{\dzC}{{C^\sharp}}
\ndef{\tD}{{\widetilde D}} \ndef{\tcD}{{\widetilde\clD}} \ndef{\ttcD}{\widetilde{\tcD}} \ndef{\sfD}{{\textsf D}} \ndef{\ttD}{\widetilde{\tD}} \ndef{\dzD}{{D^\sharp}}
\ndef{\tE}{{\widetilde E}} \ndef{\tcE}{{\widetilde\clE}} \ndef{\ttcE}{\widetilde{\tcE}} \ndef{\sfE}{{\textsf E}} \ndef{\ttE}{\widetilde{\tE}} \ndef{\dzE}{{E^\sharp}}
\ndef{\tF}{{\widetilde F}} \ndef{\tcF}{{\widetilde\clF}} \ndef{\ttcF}{\widetilde{\tcF}} \ndef{\sfF}{{\textsf F}} \ndef{\ttF}{\widetilde{\tF}} \ndef{\dzF}{{F^\sharp}}
\ndef{\tG}{{\widetilde G}} \ndef{\tcG}{{\widetilde\clG}} \ndef{\ttcG}{\widetilde{\tcG}} \ndef{\sfG}{{\textsf G}} \ndef{\ttG}{\widetilde{\tG}} \ndef{\dzG}{{G^\sharp}}
\ndef{\tH}{{\widetilde H}} \ndef{\tcH}{{\widetilde\clH}} \ndef{\ttcH}{\widetilde{\tcH}} \ndef{\sfH}{{\textsf H}} \ndef{\ttH}{\widetilde{\tH}} \ndef{\dzH}{{H^\sharp}}
\ndef{\tI}{{\widetilde I}} \ndef{\tcI}{{\widetilde\clI}} \ndef{\ttcI}{\widetilde{\tcI}} \ndef{\sfI}{{\textsf I}} \ndef{\ttI}{\widetilde{\tI}} \ndef{\dzI}{{I^\sharp}}
\ndef{\tJ}{{\widetilde J}} \ndef{\tcJ}{{\widetilde\clJ}} \ndef{\ttcJ}{\widetilde{\tcJ}} \ndef{\sfJ}{{\textsf J}} \ndef{\ttJ}{\widetilde{\tJ}} \ndef{\dzJ}{{J^\sharp}}
\ndef{\tK}{{\widetilde K}} \ndef{\tcK}{{\widetilde\clK}} \ndef{\ttcK}{\widetilde{\tcK}} \ndef{\sfK}{{\textsf K}} \ndef{\ttK}{\widetilde{\tK}} \ndef{\dzK}{{K^\sharp}}
\ndef{\tL}{{\widetilde L}} \ndef{\tcL}{{\widetilde\clL}} \ndef{\ttcL}{\widetilde{\tcL}} \ndef{\sfL}{{\textsf L}} \ndef{\ttL}{\widetilde{\tL}} \ndef{\dzL}{{L^\sharp}}
\ndef{\tM}{{\widetilde M}} \ndef{\tcM}{{\widetilde\clM}} \ndef{\ttcM}{\widetilde{\tcM}} \ndef{\sfM}{{\textsf M}} \ndef{\ttM}{\widetilde{\tM}} \ndef{\dzM}{{M^\sharp}}
\ndef{\tN}{{\widetilde N}} \ndef{\tcN}{{\widetilde\clN}} \ndef{\ttcN}{\widetilde{\tcN}} \ndef{\sfN}{{\textsf N}} \ndef{\ttN}{\widetilde{\tN}} \ndef{\dzN}{{N^\sharp}}
\ndef{\tO}{{\widetilde O}} \ndef{\tcO}{{\widetilde\clO}} \ndef{\ttcO}{\widetilde{\tcO}} \ndef{\sfO}{{\textsf O}} \ndef{\ttO}{\widetilde{\tO}} \ndef{\dzO}{{O^\sharp}}
\ndef{\tP}{{\widetilde P}} \ndef{\tcP}{{\widetilde\clP}} \ndef{\ttcP}{\widetilde{\tcP}} \ndef{\sfP}{{\textsf P}} \ndef{\ttP}{\widetilde{\tP}} \ndef{\dzP}{{P^\sharp}}
\ndef{\tQ}{{\widetilde Q}} \ndef{\tcQ}{{\widetilde\clQ}} \ndef{\ttcQ}{\widetilde{\tcQ}} \ndef{\sfQ}{{\textsf Q}} \ndef{\ttQ}{\widetilde{\tQ}} \ndef{\dzQ}{{Q^\sharp}}
\ndef{\tR}{{\widetilde R}} \ndef{\tcR}{{\widetilde\clR}} \ndef{\ttcR}{\widetilde{\tcR}} \ndef{\sfR}{{\textsf R}} \ndef{\ttR}{\widetilde{\tR}} \ndef{\dzR}{{R^\sharp}}
\ndef{\tS}{{\widetilde S}} \ndef{\tcS}{{\widetilde\clS}} \ndef{\ttcS}{\widetilde{\tcS}} \ndef{\sfS}{{\textsf S}} \ndef{\ttS}{\widetilde{\tS}} \ndef{\dzS}{{S^\sharp}}
\ndef{\tT}{{\widetilde T}} \ndef{\tcT}{{\widetilde\clT}} \ndef{\ttcT}{\widetilde{\tcT}} \ndef{\sfT}{{\textsf T}} \ndef{\ttT}{\widetilde{\tT}} \ndef{\dzT}{{T^\sharp}}
\ndef{\tU}{{\widetilde U}} \ndef{\tcU}{{\widetilde\clU}} \ndef{\ttcU}{\widetilde{\tcU}} \ndef{\sfU}{{\textsf U}} \ndef{\ttU}{\widetilde{\tU}} \ndef{\dzU}{{U^\sharp}}
\ndef{\tV}{{\widetilde V}} \ndef{\tcV}{{\widetilde\clV}} \ndef{\ttcV}{\widetilde{\tcV}} \ndef{\sfV}{{\textsf V}} \ndef{\ttV}{\widetilde{\tV}} \ndef{\dzV}{{V^\sharp}}
\ndef{\tW}{{\widetilde W}} \ndef{\tcW}{{\widetilde\clW}} \ndef{\ttcW}{\widetilde{\tcW}} \ndef{\sfW}{{\textsf W}} \ndef{\ttW}{\widetilde{\tW}} \ndef{\dzW}{{W^\sharp}}
\ndef{\tX}{{\widetilde X}} \ndef{\tcX}{{\widetilde\clX}} \ndef{\ttcX}{\widetilde{\tcX}} \ndef{\sfX}{{\textsf X}} \ndef{\ttX}{\widetilde{\tX}} \ndef{\dzX}{{X^\sharp}}
\ndef{\tY}{{\widetilde Y}} \ndef{\tcY}{{\widetilde\clY}} \ndef{\ttcY}{\widetilde{\tcY}} \ndef{\sfY}{{\textsf Y}} \ndef{\ttY}{\widetilde{\tY}} \ndef{\dzY}{{Y^\sharp}}
\ndef{\tZ}{{\widetilde Z}} \ndef{\tcZ}{{\widetilde\clZ}} \ndef{\ttcZ}{\widetilde{\tcZ}} \ndef{\sfZ}{{\textsf Z}} \ndef{\ttZ}{\widetilde{\tZ}} \ndef{\dzZ}{{Z^\sharp}}

\ndef{\bfa}{{\mathbf a}}
\ndef{\bfb}{{\mathbf b}}
\ndef{\bfc}{{\mathbf c}}
\ndef{\bfd}{{\mathbf d}}

\ndef{\euu}{{\EuScript u}}

  \ndef{\eps}{\varepsilon}

\let\geq\geqslant

\ndef{\lims}[1]{\lim\limits_{#1}}
\ndef{\sums}[1]{\sum\limits_{#1}}
\ndef{\ints}[1]{\int_{#1}}
\ndef{\sups}[1]{\sup\limits_{#1}}
\ndef{\liminfty}[1]{\lims{#1\to\infty}}
\ndef{\suminf}[1]{\sums{#1=1}^\infty}

\ndef{\limo}[1]{\omega\mbox{-}\!\!\!\lims{#1\to\infty}}          
\ndef{\limL}[1]{\rmL\mbox{-}\!\!\!\lims{#1\to\infty}}            
\ndef{\limLOne}[1]{\clL_1\mbox{-}\!\lims{#1}}
\ndef{\tildelimo}[1]{\tilde\omega\mbox{-}\!\!\!\lims{#1\to\infty}}
\ndef{\slim}{\mathrm{s}\mbox{-}\!\!\lim}          
\ndef{\wlim}{\mathrm{w}\mbox{-}\!\lim}          

\ndef{\Aut}{\operatorname{Aut}}      
\ndef{\Ch}{\operatorname{ch}}        
\ndef{\End}{\operatorname{End}}      
\ndef{\Hom}{\operatorname{Hom}}      
\rndef{\ker}{\operatorname{ker}}      
\ndef{\coker}{\operatorname{coker}}      
\ndef{\im}{\operatorname{im}}        
\ndef{\Log}{\operatorname{Log}}      
\ndef{\OP}{\operatorname{OP}}        
\ndef{\Op}{\operatorname{Op}}        
\ndef{\Symb}{\operatorname{Symb}}    
\ndef{\Wres}{\operatorname{Wres}}    
\ndef{\cl}{\operatorname{cl}}        
\ndef{\com}{\operatorname{com}}
\ndef{\const}{\operatorname{const}}  
\ndef{\conv}{\operatorname{conv}}    
\ndef{\Var}{\operatorname{Var}}
\ndef{\Cov}{\operatorname{Cov}}

\ndef{\detFK}[1]{\Delta\brs{#1}} 
\ndef{\detFKrel}[2]{\Delta_{#2}\brs{#1}} 

\ndef{\adj}{\operatorname{adj}}    
\ndef{\diag}{\operatorname{diag}}    
\ndef{\dist}{\operatorname{dist}}    
\ndef{\dom}{\operatorname{dom}}      
\ndef{\ec}{\operatorname{ec}}        
\ndef{\id}{\mathrm{Id}}                        
\ndef{\ind}{\operatorname{ind}}      
\ndef{\mydeg}{\operatorname{deg}}    
\ndef{\op}{\operatorname{op}}
\ndef{\rank}{\operatorname{rank}}
\ndef{\res}{\operatorname{res}}      
\ndef{\ran}{\operatorname{ran}}      
\ndef{\sflow}{\operatorname{sf}}     
\ndef{\isf}{\operatorname{isf}}      
\ndef{\sign}{\operatorname{sign}}    
\ndef{\sgn}{\operatorname{sgn}}      
\ndef{\sing}{\operatorname{sing}}    
\ndef{\supp}{\operatorname{supp}}    
\ndef{\tr}{\operatorname{tr}}        
\ndef{\var}{\operatorname{var}}      
\ndef{\vol}{\operatorname{vol}}      
\ndef{\wn}{\operatorname{wn}}        
\ndef{\wres}{\operatorname{wres}}    

\ndef{\prng}[1]{\mathrm R_{#1}} 
\ndef{\pker}[1]{\mathrm N_{#1}} 
\ndef{\rprng}[2]{\mathrm R_{#1}^{#2}}           
\ndef{\rpker}[2]{\mathrm N_{#1}^{#2}}           
\ndef{\rsupp}[1]{\supp_r(#1)}
\ndef{\lsupp}[1]{\supp_l(#1)}
\ndef{\rslv}[1]{R_z(#1)}      
\ndef{\HH}{H}                 
\ndef{\tHH}{\tilde \HH}       
\ndef{\VV}{V}                 
\ndef{\Rz}{R_z}               
\ndef{\tRz}{\tR_z}            
\ndef{\psif}[1]{#1^{[1]}} 
\ndef{\WPlus}[1]{W_{#1}(\mbR)} 

\ndef{\bndl}{\xi}                         
\ndef{\bndlA}{\eta}                       
\ndef{\GlueMap}{\varphi}                  
\ndef{\ChartMap}{h}                       
\ndef{\chern}{\ensuremath{\mathrm{ch}}}
\ndef{\hilba}{\clH^{(a)}}                    
\ndef{\hilbs}{\clH^{(s)}}                    
   \ndef{\hilbasargument}{(\hilb)} 
\ndef{\LpH}[1]{\clL_{#1}\hilbasargument}       
\ndef{\saLpH}[1]{\clL_{sa}^{#1}\hilbasargument}       
\ndef{\clBH}{\clB\hilbasargument}              
\ndef{\ubBH}{\clB_1\hilbasargument}            
\ndef{\clCH}{\clC\hilbasargument}              
\ndef{\clKH}{\clK\hilbasargument}              
\ndef{\clFH}{\clF\hilbasargument}              
\ndef{\clUH}{\clU\hilbasargument}              
\ndef{\clCFH}{{\clC\clF}\hilbasargument}       
\ndef{\saBH}{\clB_{sa}\hilbasargument}         
\ndef{\saCH}{\clC_{sa}\hilbasargument}         
\ndef{\saFH}{\clF_{sa}\hilbasargument}         
\ndef{\saKH}{\clK_{sa}\hilbasargument}         
\ndef{\saCFH}{\clC\clF_{sa}\hilbasargument}    
\ndef{\clUFH}{\clU\clF\hilbasargument}         
\ndef{\Uinj}{\clU_{inj}\hilbasargument}        
\ndef{\UFinj}{\clU\clF_{inj}\hilbasargument}   

\ndef{\spproj}[2]{E^{#1}_{#2}}                      
\ndef{\spprojb}[2]{E^{#2}_{#1}}                     

\ndef{\LpN}[1]{\clL^{#1}(\clN,\tau)}     
\ndef{\saLpN}[1]{\clL^{#1}_{sa}(\clN,\tau)} 
\ndef{\rLpN}[1]{L^{#1}(\clN,\tau)}       
\ndef{\clAND}{(\clA,\clN,D)}             
\ndef{\clBA}{{\clB(\clA)}}
\ndef{\saKN}{{\clK_{sa}(\clN,\tau)}}          
\ndef{\clKN}{{\clK(\clN,\tau)}}          
\ndef{\clKtN}{{\clK(\tilde\clN,\tau)}}   
\ndef{\clFN}{{\clF(\clN,\tau)}}          
\ndef{\saFN}{{\clF_{sa}(\clN,\tau)}}     
\ndef{\clPN}{\clP(\clN)}                 
\ndef{\clQN}{\clQ(\clN,\tau)}            
\ndef{\infPN}{{\clP_\tau^\infty(\clN)}}  
\ndef{\clOF}[2]{\clF_{#1\mbox{-}#2}(\clN,\tau)}         
\ndef{\oind}[2]{{\rm \tau\mbox{-}ind}_{#1\mbox{-}#2}}   
\ndef{\tind}{\tau\mbox{-}\ind}                  
\ndef{\DInd}{\ind_{\clD,\tau}}           
\ndef{\BF}{Breuer-Fredholm}              
\ndef{\skewfred}[2]{$(#1\cdot #2)$ $\tau$\tire Fredholm}   
\ndef{\affl}{\eta}                       
\ndef{\vNa}{von Neumann algebra}         
\ndef{\nsf}{faithful normal semifinite } 
\ndef{\taubrs}[1]{\tau\brackets{#1}}     
\ndef{\sqbrs}[1]{[#1]}        
\ndef{\Sqbrs}[1]{\big[#1\big]}        
\ndef{\SqBrs}[1]{\Big[#1\Big]}        

\ndef{\domd}{\bigcap\limits_{n\ge 0} \dom\;\delta^n}         
\ndef{\DiffOP}{{\rm \clD}}
\ndef{\ADA}{\clA \cup [\clD,\clA]}
\ndef{\DixIdeal}[1]{\LpH{#1,\infty}}               
\ndef{\dixideal}{\ell^{1,\infty}}                  
\ndef{\WDixIdeal}{\LpH{1,\mathrm w}}               
\ndef{\DixIdealPos}[1]{\DixIdeal{#1}_+}            
\ndef{\DixIdealN}[1]{\LpN{#1,\infty}}              
\ndef{\DixIdealNPar}[2]{\clL^{#1,\infty}_{#2}(\clN,\tau)}    
\ndef{\DixIdealNPos}[1]{\LpN{#1,\infty}_+}                   
\ndef{\TrD}{\Tr_\omega}                                      
\ndef{\tauD}{{\tau_\omega}}                                  
\ndef{\ILogN}{\frac 1{\log(1+N)}}
\ndef{\DixNorm}[1]{\norm{#1}_{(1,\infty)}}                   
\ndef{\DixInt}[1]{\ints 0^t \mu_s(#1)\,ds}
\ndef{\DixIntL}[1]{\ints 0^{\lambda_{1/t}(#1)}\mu_s(#1)\,ds}
    \ndef{\SmallIdeal}{{\clL^{1, \mathrm w}}}
    \ndef{\SmallIdealMeas}{{\clL^{1, \mathrm w}_m}}
    \ndef{\DixIntII}[1]{\int_0^t \mu_s(#1)\,ds}
    \ndef{\DixIntf}[1]{\Phi_t(#1)}
    \ndef{\DixIntg}[1]{\Psi_t(#1)}

\ndef{\lpi}{\clL^{1,\pi}(\clN,\tau)}

\ndef{\strl}[1]{\stackrel \longrightarrow {#1}}
\ndef{\IIinfty}{$\mathrm{II}_\infty$\ }

\ndef{\fourier}[1]{\clF(#1)}          
\ndef{\HaarMeasBohrs}{\nu}            
\ndef{\BrownsMeas}{\mu}               
\ndef{\BohrCont}[1]{\tilde{#1}}       
\ndef{\APMean}{{M}}                   
\ndef{\CDSS}{{\clA_B}}                
\ndef{\matr}{{\rm Mat}}               
\ndef{\seque}[1]{\ensuremath{\{#1_n\}_{n=1}^\infty}}    
\ndef{\sequen}[2]{\ensuremath{\{#1_#2\}_{#2=1}^\infty}}    
\ndef{\Seque}[1]{\ensuremath{\left(#1_0,#1_1,#1_2,\dots\right)}}    
\ndef{\Cesaro}{H}                           
\ndef{\CesaroRPlus}{M}                      
\ndef{\Dilation}{D}                         
\ndef{\Shift}{T}                            

\ndef{\TrNorm}[1]{\norm{#1}_1}              
\ndef{\HSNorm}[1]{\norm{#1}_2}              
\ndef{\InftyNorm}[1]{\norm{#1}_\infty}      
\ndef{\normQN}[1]{\norm{#1}_{\clQN}}        
\ndef{\clLpnorm}[2]{\norm{#2}_{\clL^{#1}}}    
\ndef{\clLnorm}[1]{\clLpnorm{1}{#1}}    

\ndef{\ccurve}{\gamma}                      

\ndef{\Brs}[1]{\big(#1\big)}                
\ndef{\BRS}[1]{\Big(#1\Big)}                
\ndef{\scal}[2]{\left\langle #1,#2\right\rangle}               
\ndef{\Scal}[1]{\left\langle #1\right\rangle}               
\ndef{\precprec}{\prec\!\!\!\prec}
\ndef{\qeq}{\stackrel?=}
\ndef{\spectrum}[1]{\sigma_{#1}} 
\ndef{\spectruma}[1]{\sigma^{(a)}_{#1}} 
\rndef{\emptyset}{\varnothing}                              
\ndef{\csupp}{c}                           
\ndef{\closure}[1]{\overline{#1}}
\ndef{\linspan}[1]{\mathrm{span}\,{#1}}
\ndef{\bddborel}[1]{B(#1)}                 
\ndef{\charfunc}{\chi}
\ndef{\FrDer}{\euD}                        
\ndef{\LieDer}[1]{\pounds_{#1}\,}          
\ndef{\dds}{\left.\frac d{ds} \right|_{s = 0}}
\ndef{\ortcmp}[1]{#1^{\scriptscriptstyle \perp}}            
\ndef{\Laplace}{\Delta}                    

\ndef{\matrPQ}[3]
{
    \left(
      \begin{array}{cc}
        #1_{11} & #1_{12} \\
        #1_{21} & #1_{22}
      \end{array}
    \right)_{[#2,#3]}
}

\ndef{\margOK}{\marginpar{\bf \small OK}}

\newcounter{margcomcount}
\ndef{\margcom}[1]{\marginpar{\bf \small #1} \addtocounter{margcomcount}{1}
   \index{\indexcom{{\bf COMMENT: #1}}}}

\ndef{\mytimes}{\!\times\!}
\ndef{\sss}[1]{\subsubsection{}\label{#1}}

\rndef{\phi}{\varphi} \ndef{\OpenUnitDisk}{D}
\ndef{\RHS}{RHS}                            
\ndef{\LHS}{LHS} 
\ndef{\ttt}{\Leftrightarrow}
\ndef{\then}{\Rightarrow}
\ndef{\tto}{\longrightarrow}
\ndef{\nno}{\nonumber\\}
\ndef{\newn}[1]{\index{#1} {\bfseries #1}}       
\ndef{\la}{\langle}
\ndef{\ra}{\rangle}
\ndef{\dbar}{{\;\bar{\phantom{o}} \!\!\!\! d}}
\ndef{\stl}[1]{\stackrel{\vbox to 0pt{\vss\hbox{$\scriptstyle #1$}}}}
\ndef{\mathcomment}[1]{{\hfill \qquad\qquad\qquad\text{by (#1)}}}        
\ndef{\mathcomm}[1]{{\hfill \qquad\qquad\qquad\qquad\qquad\text{#1}}}        
\ndef{\details}[1]{\smallskip\begin{center} {\bf Here:}
#1\end{center}\medskip} \ndef{\indexcom}[1]{ --- #1}
\ndef{\longsim}{\ \sim \ }              
\ndef{\tire}{-}              
\ndef{\intinfinf}{\int_{-\infty}^\infty}
     \ndef{\npartial}{\slash\!\!\!\partial}
     \ndef{\Heis}{\operatorname{Heis}}
     \ndef{\Solv}{\operatorname{Solv}}
     \ndef{\Spin}{\operatorname{Spin}}
     \ndef{\SO}{\operatorname{SO}}
     \ndef{\Index}{\operatorname{index}}

             \ndef{\p}{\partial}
             \ndef{\dd}{|\clD|}
             \ndef{\n}{\parallel}

\let\LatexCite=\cite  

\let\ifnumref\iffalse 

\ndef{\ifuncited}[4]{\expandafter\ifx\csname used#4\endcsname\relax}

\ndef{\ifcited}[4]{\expandafter\ifx\csname used#4\endcsname\relax\else}

  \ndef{\papertitle}[1]{ \emph{#1}, }
  \ndef{\paperauthor}[2]{#2}  
  \ndef{\pbbi}[9]{%
      \ifcited{#1}{#2}{#3}{#5}%
        \ifnumref%
          \bibitem{#5}\paperauthor{#1}{#6},\papertitle{#7}#8.%
        \else%
          \advance #9 by 1%
          \ifnum#9<1%
            \bibitem[#4]{#5}\paperauthor{#1}{#6}, \papertitle{#7}#8.%
          \else%
            \bibitem[#4{\the#9}]{#5}\paperauthor{#1}{#6},\papertitle{#7}#8.%
          \fi%
        \fi%
      \fi%
  }
  \ndef{\mbbi}[8]{%
     \ifcited{#1}{#2}{#3}{#5}%
        \ifnumref%
          \bibitem{#5}\paperauthor{#1}{#6},\papertitle{#7}#8.%
        \else%
          \bibitem[#4]{#5}\paperauthor{#1}{#6},\papertitle{#7}#8.%
        \fi%
     \fi%
  }

\ndef{\AddCite}[1]{%
   \ifuncited{0}{0}{0}{#1}%
     \expandafter\gdef\csname used#1\endcsname {}%
   \fi%
}

\def\ProcessCite#1,{%
     \ifx\relax#1%
         \let\next=\relax%
     \else%
         \AddCite{#1}%
         \let\next=\ProcessCite%
     \fi%
     \next%
}

\ndef{\AddCites}[1]{\ProcessCite#1,\relax,}

\ndef{\CiteWithoutExtension}[1]{%
   \AddCites{#1}%
   \LatexCite{#1}%
}

\def\CiteWithExtension[#1]#2{%
   \AddCites{#2}%
   \LatexCite[#1]{#2}%
}

\ndef{\CleverCite}{%
    \ifx\NChar[ %
       \let\MyCite=\CiteWithExtension %
    \else %
       \let\MyCite=\CiteWithoutExtension %
    \fi %
    \MyCite%
}

\renewcommand{\cite}{\futurelet\NChar\CleverCite}

      \ndef{\volume}[1]{{\bf #1}}
      \ndef{\VolYearPP}[3]{\ifnum#2=0 (to appear)\else\volume{#1} (#2), #3\fi}
      \ndef{\VolNoYearPP}[4]{\ifnum#3=0 (to appear)\else\volume{#1} #2 (#3), #4\fi}
      \ndef{\libcode}[1]{}

\ndef{\jnActaMath}[3]{Acta Math. \VolYearPP{#1}{#2}{#3}}                       
\ndef{\jnAdvMath}[3]{Adv.\,in~Math. \VolYearPP{#1}{#2}{#3}}                     
\ndef{\jnAlgAnal}[3]{Algebra i~Analiz \VolYearPP{#1}{#2}{#3}}
\ndef{\jnAmerJMath}[3]{Amer.\,J.\,Math. \VolYearPP{#1}{#2}{#3}}                  
\ndef{\jnAmerMathMonth}[3]{Amer.\,Math.\,Monthly \VolYearPP{#1}{#2}{#3}}         
\ndef{\jnAnnMath}[4]{Ann. of~Math. \VolNoYearPP{#1}{#2}{#3}{#4}}               
\ndef{\jnAnalMath}[3]{J. Anal. Math. \VolYearPP{#1}{#2}{#3}}                   
\ndef{\jnArchRatMechAnal}[3]{Arch. Rational Mech. Anal. \VolYearPP{#1}{#2}{#3}}                   
\ndef{\jnBullLondMathSoc}[3]{Bull. London Math. Soc. \VolYearPP{#1}{#2}{#3}}   
\ndef{\jnBullAMS}[3]{Bull. Amer. Math. Soc. \VolYearPP{#1}{#2}{#3}}   
\ndef{\jnCanMathBull}[3]{Canad. Math. Bull. \VolYearPP{#1}{#2}{#3}}            
\ndef{\jnCanMath}[3]{Canad. J.~Math. \VolYearPP{#1}{#2}{#3}}             
\ndef{\jnCommMathPhys}[3]{Comm. Math. Phys. \VolYearPP{#1}{#2}{#3}}             
\ndef{\jnCommPDE}[3]{Comm. Partial Differential Equations \VolYearPP{#1}{#2}{#3}}             
\ndef{\jnComptRendue}[3]{C.\,R.~Acad. Sci. Paris S\'er. A-B \VolYearPP{#1}{#2}{#3}}      
\ndef{\jnContMath}[3]{Contemporary Math. \VolYearPP{#1}{#2}{#3}}               %
\ndef{\jnDukeMJ}[3]{Duke Math. J. \VolYearPP{#1}{#2}{#3}}
\ndef{\jnDiffGeom}[3]{J.~Diff. Geom. \VolYearPP{#1}{#2}{#3}}                   
\ndef{\jnErgodicTheory}[3]{Ergodic Theory and Dynamical Systems \VolYearPP{#1}{#2}{#3}} 
\ndef{\jnFuncAnal}[3]{J.~Functional Analysis \VolYearPP{#1}{#2}{#3}}           
\ndef{\jnFunkAnalPril}[4]{Funct. Anal. Appl. \VolNoYearPP{#1}{#2}{#3}{#4}}  
\ndef{\jnGAFA}[3]{GAFA \VolYearPP{#1}{#2}{#3}}                                 
\ndef{\jnIHES}[3]{IHES Publ. Math. (Paris) \VolYearPP{#1}{#2}{#3}}             
\ndef{\jnIEOT}[3]{Integral Equations Operator Theory   \VolYearPP{#1}{#2}{#3}} 
\ndef{\jnIsrMath}[3]{Israel J.~Math. \VolYearPP{#1}{#2}{#3}}                   
\ndef{\jnKTheory}[3]{K-Theory \VolYearPP{#1}{#2}{#3}}                          
\ndef{\jnLetMathPhys}[3]{Lett. Math. Phys. \VolYearPP{#1}{#2}{#3}}             
\ndef{\jnMathAnn}[3]{Math. Ann. \VolYearPP{#1}{#2}{#3}}                        
\ndef{\jnMathAnalAppl}[3]{J.~Math. Anal. and Appl. \VolYearPP{#1}{#2}{#3}}     
\ndef{\jnMathNachr}[3]{Math.\,Nachr. \VolYearPP{#1}{#2}{#3}}
\ndef{\jnMathPhys}[3]{J. Math. Phys. \VolYearPP{#1}{#2}{#3}}
\ndef{\jnMathSocJap}[3]{J. Math. Soc. Japan \VolYearPP{#1}{#2}{#3}}
\ndef{\jnOperTheory}[3]{J.~Operator Theory \VolYearPP{#1}{#2}{#3}}             
\ndef{\jnPacJMath}[3]{Pacific J.~Math. \VolYearPP{#1}{#2}{#3}}                  
\ndef{\jnPositivity}[3]{Positivity \VolYearPP{#1}{#2}{#3}}
\ndef{\jnProcAmerMS}[3]{Proc. Amer. Math. Soc. \VolYearPP{#1}{#2}{#3}}         
\ndef{\jnProcCambPhilSoc}[3]{Math. Proc. Camb. Phil. Soc. \VolYearPP{#1}{#2}{#3}}
\ndef{\jnReineAngew}[3]{J.~Reine Angew. Math. \VolYearPP{#1}{#2}{#3}}          
\ndef{\jnTokyoMath}[3]{Tokyo J.~Math. \VolYearPP{#1}{#2}{#3}}
\ndef{\jnTopology}[3]{Topology \VolYearPP{#1}{#2}{#3}}
\ndef{\jnTransAmerMathSoc}[3]{Trans. Amer. Math. Soc. \VolYearPP{#1}{#2}{#3}}
\ndef{\jnIzvANSSSR}[3]{Izv. Akad. Nauk SSSR, Ser. Mat. \VolYearPP{#1}{#2}{#3}}
\ndef{\jnIzvVyshUchZav}[3]{Izv. Vyssh. Uch. Zav., Mat. \VolYearPP{#1}{#2}{#3} (Russian)}
\ndef{\jnIzdatLenUniv}[2]{Izdat. Leningrad. Univ., Leningrad, (#1), #2 (Russian)}
\ndef{\jnFieldsInsComm}[3]{Fields Inst. Comm. \VolYearPP{#1}{#2}{#3}}
\ndef{\jnDoklANSSSR}[3]{Dokl. Akad. Nauk SSSR \VolYearPP{#1}{#2}{#3}}
\ndef{\jnMatZametki}[3]{Matem. zametki \VolYearPP{#1}{#2}{#3}}
\ndef{\jnRussMathSurvey}[3]{Russian Math. Surveys \VolYearPP{#1}{#2}{#3}}
\ndef{\jnSibMathJ}[3]{Sib. Math.~J. \VolYearPP{#1}{#2}{#3}}
\ndef{\jnSovMath}[3]{J.~Soviet math. \VolYearPP{#1}{#2}{#3}}
\ndef{\jnTransMoscMathSoc}[3]{Trans. Moscow Math. Soc. \VolYearPP{#1}{#2}{#3}}
\ndef{\jnUMN}[3]{Uspekhi Mat. Nauk \VolYearPP{#1}{#2}{#3}}

\ndef{\bkTransMathMon}[2]{Trans. Math. Monographs, AMS, \volume{#1}, #2}

\ndef{\pbBirkhauser}[1]{Birkh\"auser, Boston, #1}
\ndef{\pbFactorial}[1]{Moscow, Factorial, #1}
\ndef{\pbGauthier}[1]{Gauthier-Villars, Paris, #1}
\ndef{\pbNauka}[1]{Moscow, Nauka, #1 (Russian)}
\ndef{\pbNaukaR}[1]{Москва, Наука, #1}
\ndef{\pbPrinceton}[1]{Princeton University Press, Princeton, New Jersey, #1}
\ndef{\pbPublPerish}[1]{Publish or Perish Inc., Berkeley, #1}
\ndef{\pbSpringer}[1]{Springer-Verlag, #1}

\ndef{\myauthor}[1]{\mbox{#1}}

\ndef{\Agmon}{\myauthor{Sh.\,Agmon}}
\ndef{\Ahiezer}{\myauthor{N.\,I.\,Ahiezer}}
\ndef{\Arazy}{\myauthor{J.\,Arazy}}
\ndef{\Aronszajn}{\myauthor{N.\,Aronszajn}}
\ndef{\Astashkin}{\myauthor{S.\,V.\,Astashkin}}
\ndef{\Atiyah}{\myauthor{M.\,Atiyah}}
\ndef{\Avron}{\myauthor{J.\,E.\,Avron}}
\ndef{\Azamov}{\myauthor{N.\,A.\,Azamov}}
\ndef{\Banach}{\myauthor{S.\,Banach}}
\ndef{\Benameur}{\myauthor{M-T.\,Benameur}}
\ndef{\Bennett}{\myauthor{C.\,Bennett}}
\ndef{\Berezin}{\myauthor{F.\,A.\,Berezin}}
\ndef{\Berline}{\myauthor{N.\,Berline}}
\ndef{\Birman}{\myauthor{M.\,Sh.\,Birman}}
\ndef{\Blackadar}{\myauthor{B.\,Blackadar}}
\ndef{\Bogolyubov}{\myauthor{N.\,N.\,Bogolyubov}}
\ndef{\Bonsall}{\myauthor{F.\,F.\,Bonsall}}
\ndef{\Bony}{\myauthor{J.\,F.\,Bony}}
\ndef{\BoosBavnbek}{\myauthor{B.\,Boo$\beta$-Bavnbek}}
\ndef{\Bott}{\myauthor{R.\,Bott}}
\ndef{\Branges}{\myauthor{L.\,de Branges}}
\ndef{\Bratteli}{\myauthor{O.\,Bratteli}}
\ndef{\Bredon}{\myauthor{G.\,E.\,Bredon}}
\ndef{\Breuer}{\myauthor{M.\,Breuer}}
\ndef{\Brown}{\myauthor{L.\,G.\,Brown}}
\ndef{\Bruneau}{\myauthor{V.\,Bruneau}}
\ndef{\Buslaev}{\myauthor{V.\,S.\,Buslaev}}
\ndef{\Carey}{\myauthor{A.\,L.\,Carey}}
\ndef{\CareyRW}{\myauthor{R.\,W.\,Carey}} 
\ndef{\Cartan}{\myauthor{H.\,Cartan}}
\ndef{\Chilin}{\myauthor{V.\,I.\,Chilin}}
\ndef{\Coburn}{\myauthor{L.\,A.\,Coburn}}
\ndef{\Connes}{\myauthor{A.\,Connes}}
\ndef{\Cornfeld}{\myauthor{I.\,P.\,Cornfeld}}
\ndef{\Daletskii}{\myauthor{Yu.\,L.\,Daletski\u\i}}   
\ndef{\Dixmier}{\myauthor{J.\,Dixmier}}
\ndef{\DoddsPG}{\myauthor{P.\,G.\,Dodds}}
\ndef{\DoddsTK}{\myauthor{T.\,K.\,Dodds}}
\ndef{\Douglas}{\myauthor{R.\,G.\,Douglas}}
\ndef{\Dubrovin}{\myauthor{B.\,A.\,Dubrovin}}
\ndef{\Dugundji}{\myauthor{J.\,Dugundji}}
\ndef{\Duncan}{\myauthor{J.\,Duncan}}
\ndef{\Dunford}{\myauthor{N.\,Dunford}}
\ndef{\Dykema}{\myauthor{K.\,J.\,Dykema}}
\ndef{\Edwards}{\myauthor{R.\,E.\,Edwards}}
\ndef{\Eilenberg}{\myauthor{S.\,Eilenberg}}
\ndef{\Entina}{\myauthor{S.\,B.\,\`Entina}}
\ndef{\Fack}{\myauthor{T.\,Fack}} 
\ndef{\Faddeev}{\myauthor{L.\,D.\,Faddeev}}
\ndef{\Farber}{\myauthor{M.\,Farber}}
\ndef{\Farforovskaya}{\myauthor{Yu.\,B.\,Farforovskaya}}
\ndef{\Federer}{\myauthor{H.\,Federer}}
\ndef{\Fedosov}{\myauthor{B.\,V.\,Fedosov}}
\ndef{\Figiel}{\myauthor{T.\,Figiel}} 
\ndef{\Figueroa}{\myauthor{H.\,Figueroa}}
\ndef{\Fillmore}{\myauthor{P.\,A.\,Fillmore}}
\ndef{\Fomenko}{\myauthor{A.\,T.\,Fomenko}} 
\ndef{\Fomin}{\myauthor{S.\,V.\,Fomin}}
\ndef{\Frohlich}{\myauthor{J.\,Fr\"ohlich}}
\ndef{\Fuglede}{\myauthor{B.\,Fuglede}}
\ndef{\Furutani}{\myauthor{K.\,Furutani}}
\ndef{\Gelfand}{\myauthor{I.\,M.\,Gelfand}}
\ndef{\Gesztesy}{\myauthor{F.\,Gesztesy}}     
\ndef{\Getzler}{\myauthor{E.\,Getzler}} 
\ndef{\Gilkey}{\myauthor{P.\,B.\,Gilkey}}
\ndef{\Gitler}{\myauthor{S.\,Gitler}}
\ndef{\Glazman}{\myauthor{I.\,M.\,Glazman}}
\ndef{\Glimm}{\myauthor{J.\,Glimm}}
\ndef{\Gohberg}{\myauthor{I.\,C.\,Gohberg}}
\ndef{\Goldshtein}{\myauthor{Ya.\,Goldshtein}}
\ndef{\Golze}{\myauthor{F.\,Golze}}
\ndef{\GraciaBondia}{\myauthor{J.\,M.\,Gracia-Bond\'{i}a}}
\ndef{\Greenleaf}{\myauthor{F.\,P.\,Greenleaf}}
\ndef{\Gromov}{\myauthor{M.\,Gromov}}
\ndef{\Gunning}{\myauthor{R.\,C.\,Gunning}}
\ndef{\Haagerup}{\myauthor{U.\,Haagerup}}
\ndef{\Haag}{\myauthor{R.\,Haag}}
\ndef{\Halmos}{\myauthor{P.\,R.\,Halmos}}
\ndef{\Hardy}{\myauthor{G.\,H.\,Hardy}}
\ndef{\Herbst}{\myauthor{I.\,W.\,Herbst}}
\ndef{\Higson}{\myauthor{N.\,Higson}}  
\ndef{\Hoermander}{\myauthor{L.\,H\"ormander}} 
\ndef{\Hoffman}{\myauthor{K.\,Hoffman}} 
\ndef{\Ito}{\myauthor{K.\,Ito}}
\ndef{\Ikebe}{\myauthor{T.\,Ikebe}}
\ndef{\Jaffe}{\myauthor{A.\,Jaffe}}
\ndef{\James}{\myauthor{I.\,M.\,James}}
\ndef{\Javrjan}{\myauthor{V.\,A.\,Javrjan}}
\ndef{\Jitomirskaya}{\myauthor{S.\,Jitomirskaya}}
\ndef{\Kadison}{\myauthor{R.\,V.\,Kadison}}
\ndef{\Kalton}{\myauthor{N.\,J.\,Kalton}} 
\ndef{\Kato}{\myauthor{T.\,Kato}} 
\ndef{\Kobayashi}{\myauthor{S.\,Kobayashi}}
\ndef{\Koplienko}{\myauthor{L.\,S.\,Koplienko}}
\ndef{\Korotyaev}{\myauthor{E.\,Korotyaev}}
\ndef{\Kosaki}{\myauthor{H.\,Kosaki}}
\ndef{\Kostrykin}{\myauthor{V.\,Kostrykin}}
\ndef{\Kotani}{\myauthor{S.\,Kotani}}
\ndef{\Krein}{\myauthor{Kre\u\i n}}
\ndef{\KreinMG}{\myauthor{M.\,G.\,Kre\u\i n}}
\ndef{\KreinSG}{\myauthor{S.\,G.\,Kre\u\i n}}
\ndef{\Kuroda}{\myauthor{S.\,T.\,Kuroda}}
\ndef{\Leichtnam}{\myauthor{E.\,Leichtnam}}
\ndef{\Lesch}{\myauthor{M.\,Lesch}}
\ndef{\Lesniewski}{\myauthor{A.\,Lesniewski}}
\ndef{\Levitan}{\myauthor{B.\,M.\,Levitan}}
\ndef{\Lidskii}{\myauthor{V.\,B.\,Lidskii}}
\ndef{\Lifshitz}{\myauthor{I.\,M.\,Lifshitz}}
\ndef{\Lindenstrauss}{\myauthor{J.\,Lindenstrauss}}
\ndef{\Loday}{\myauthor{J.-L.\,Loday}}
\ndef{\Lord}{\myauthor{S.\,Lord}}      
\ndef{\Lorentz}{\myauthor{G.\,Lorentz}}
\ndef{\Magnus}{\myauthor{W.\,Magnus}}
\ndef{\Makarov}{\myauthor{K.\,A.\,Makarov}}
\ndef{\MakarovN}{\myauthor{N.\,Makarov}}
\ndef{\Mathai}{\myauthor{V.\,Mathai}}         
\ndef{\McKean}{\myauthor{H.\,P.\,McKean}}
\ndef{\Mishchenko}{\myauthor{A.\,S.\,Mishchenko}}
\ndef{\Molchanov}{\myauthor{S.\,A.\,Molchanov}}
\ndef{\Moore}{\myauthor{C.\,C.\,Moore}}
\ndef{\Moscovici}{\myauthor{H.\,Moscovici}}  
\ndef{\Motovilov}{\myauthor{A.\,K.\,Motovilov}}
\ndef{\Moyer}{\myauthor{R.\,D.\,Moyer}}
\ndef{\Naboko}{\myauthor{S.\,N.\,Naboko}}
\ndef{\Narasimhan}{\myauthor{R.\,Narasimhan}}
\ndef{\Nomizu}{\myauthor{K.\,Nomizu}}
\ndef{\Novikov}{\myauthor{S.\,P.\,Novikov}}
\ndef{\Osterwalder}{\myauthor{K.\,Osterwalder}}
\ndef{\Patodi}{\myauthor{V.\,Patodi}}
\ndef{\Pagter}{\myauthor{B.\,de~Pagter}}  
\ndef{\Pastur}{\myauthor{L.\,A.\,Pastur}}  
\ndef{\Pavlov}{\myauthor{B.\,S.\,Pavlov}}
\ndef{\Pedersen}{\myauthor{G.\,K.\,Pedersen}}
\ndef{\Peller}{\myauthor{V.\,V.\,Peller}}
\ndef{\Perera}{\myauthor{V.\,S.\,Perera}}
\ndef{\Petunin}{\myauthor{Ju.\,I.\,Petunin}}
\ndef{\Phillips}{\myauthor{J.\,Phillips}}  
\ndef{\Piazza}{\myauthor{P.\,Piazza}}   
\ndef{\Pincus}{\myauthor{J.\,D.\,Pincus}}   
\ndef{\Poincare}{Poincar\'e}
\ndef{\Postnikov}{\myauthor{M.\,M.\,Postnikov}} 
\ndef{\Povzner}{\myauthor{A.\,Ya.\,Povzner}}
\ndef{\Prinzis}{\myauthor{R.\,Prinzis}}
\ndef{\Privalov}{\myauthor{I.\,I.\,Privalov}}
\ndef{\Pushnitski}{\myauthor{A.\,B.\,Pushnitski}} 
\ndef{\Raeburn}{\myauthor{I.\,Raeburn}}
\ndef{\Raikov}{\myauthor{G.\,Raikov}}
\ndef{\Reed}{\myauthor{M.\,Reed}}
\ndef{\Rennie}{\myauthor{A.\,Rennie}}
\ndef{\Rickart}{\myauthor{C.\,E.\,Rickart}}
\ndef{\Riesz}{\myauthor{F.\,Riesz}}
\ndef{\Ringrose}{\myauthor{J.\,Ringrose}}
\ndef{\Rio}{\myauthor{R.\,del Rio}}
\ndef{\Robinson}{\myauthor{D.\,Robinson}}
\ndef{\Rossi}{\myauthor{H.\,Rossi}}
\ndef{\Rudin}{\myauthor{W.\,Rudin}}
\ndef{\Ruelle}{\myauthor{D.\,Ruelle}}
\ndef{\Ruzhansky}{\myauthor{M.\,Ruzhansky}}
\ndef{\Sakai}{\myauthor{Sh.\,Sakai}}
\ndef{\Sargsjan}{\myauthor{I.\,S.\,Sargsjan}}
\ndef{\Sato}{\myauthor{H.\,Sato}}
\ndef{\Schaeffer}{\myauthor{D.\,G.\,Schaeffer}}
\ndef{\Schluchtermann}{\myauthor{G.\,Schluchtermann}}
\ndef{\Schochet}{\myauthor{C.\,Schochet}}
\ndef{\SchroedingerE}{\myauthor{E.\,Schr\"odinger}}
\ndef{\Schroedinger}{\myauthor{Schr\"odinger}}
\ndef{\Schrohe}{\myauthor{E.\,Schrohe}}
\ndef{\Schwartz}{\myauthor{J.\,T.\,Schwartz}}
\ndef{\Sedaev}{\myauthor{A.\,A.\,Sedaev}}
\ndef{\Seiler}{\myauthor{R.\,Seiler}}
\ndef{\Semenov}{\myauthor{E.\,M.\,Semenov}}
\ndef{\Shabat}{\myauthor{B.\,V.\,Shabat}}
\ndef{\Shafarevich}{\myauthor{I.\,R.\,Shafarevich}}
\ndef{\Sharpley}{\myauthor{R.\,Sharpley}}
\ndef{\Shilov}{\myauthor{G.\,E.\,Shilov}}
\ndef{\Shirkov}{\myauthor{D.\,V.\,Shirkov}}
\ndef{\Shubin}{\myauthor{M.\,A.\,Shubin}}
\ndef{\Silverman}{\myauthor{H.\,Silverman}}
\ndef{\Simon}{\myauthor{B.\,Simon}}
\ndef{\Sinai}{\myauthor{Ya.\,G.\,Sinai}}
\ndef{\Singer}{\myauthor{I.\,M.\,Singer}}
\ndef{\Solomyak}{\myauthor{M.\,Z.\,Solomyak}}
\ndef{\Soloviev}{\myauthor{Yu.\,P.\,Soloviev}}
\ndef{\Spivak}{\myauthor{M.\,Spivak}}
\ndef{\Stein}{\myauthor{E.\,M.\,Stein}}
\ndef{\Stenkin}{\myauthor{V.\,V.\,Sten'kin}}
\ndef{\Stratila}{\myauthor{S.\,Stratila}}
\ndef{\Sucheston}{\myauthor{L.\,Sucheston}}
\ndef{\Sukochev}{\myauthor{F.\,A.\,Sukochev}}
\ndef{\Switzer}{\myauthor{R.\,M.\,Switzer}}
\ndef{\SzNagy}{\myauthor{B.\,Sz.-Nagy}}
\ndef{\Takesaki}{\myauthor{M.\,Takesaki}}
\ndef{\Taylor}{\myauthor{M.\,E.\,Taylor}}
\ndef{\Treves}{\myauthor{F.\,Treves}}
\ndef{\Troitsky}{\myauthor{E.\,V.\,Troitsky}}
\ndef{\Tzafriri}{\myauthor{L.\,Tzafriri}}
\ndef{\Varilly}{\myauthor{J.\,C.\,V\'{a}rilly}}
\ndef{\Vergne}{\myauthor{M.\,Vergne}}
\ndef{\Vladimirov}{\myauthor{V.\,S.\,Vladimirov}}
\ndef{\Voiculescu}{\myauthor{D.\,Voiculescu}}
\ndef{\Weiss}{\myauthor{G.\,Weiss}}
\ndef{\Wells}{\myauthor{R.\,O.\,Wells}}
\ndef{\Williams}{\myauthor{J.\,P.\,Williams}}
\ndef{\Winkler}{\myauthor{S.\,Winkler}}
\ndef{\Witten}{\myauthor{E.\,Witten}}
\ndef{\Wodzicki}{\myauthor{M.\,Wodzicki}}
\ndef{\Wojciechowski}{\myauthor{K.\,P.\,Wojciechowski}}
\ndef{\Yafaev}{\myauthor{D.\,R.\,Yafaev}}
\ndef{\Yosida}{\myauthor{K.\,Yosida}}
\ndef{\Zsido}{\myauthor{L.\,Zsido}}

\rndef{\iff}{\Leftrightarrow}

\newcommand{\tlV}{\tilde V}
\newcommand{\tlJ}{\tilde J}
\renewcommand{\TrD}{\mathrm{Tr}_\omega}

\begin{document}

\title[Resonance set and its structure]{Spectral flow inside essential spectrum II: \\ resonance set and its structure}

\author{Nurulla Azamov}

\address{Independent scholar, Adelaide, SA, Australia}

\email{azamovnurulla@gmail.com}
 \keywords{Spectral flow, essential spectrum}
 \subjclass[2000]{ 
     Primary 47A40}

\begin{abstract}
This paper is a continuation of the study of spectral flow inside essential spectrum initiated in \cite{AzSFIES}.
Given a point $\lambda$ outside the essential spectrum of a self-adjoint operator $H_0,$ the resonance set, $\euR(\lambda),$ 
is an analytic variety which 
consists of self-adjoint relatively compact perturbations $H_0+V$ of $H_0,$ for which $\lambda$ is an eigenvalue. 
One may ask for criteria for the vector $V$ to be tangent to the resonance set. Such criteria were given in \cite{AzSFnRI}.

In this paper we study similar criteria for the case of $\lambda$ inside the essential spectrum of $H_0.$
For the case $\lambda \in \sigma_{ess}(H_0)$ the resonance set is defined in terms of the well-known limiting absorption principle.
Among the results of this paper is that the resonance set contains plenty of straight lines, moreover, given any regular relatively compact perturbation $V$
there exists a finite rank self-adjoint operator, $\tilde V,$ such that the straight line $H_0 + \mbR(V-\tilde V)$ belongs to the resonance set. 
 
 Another result of this paper is that inside the essential spectrum there exist plenty of transversal to the resonance set perturbations $V$ which have order $\geq 2,$
 in contrast to what happens outside the essential spectrum, \cite{AzSFnRI}.


\end{abstract}

\maketitle


\setcounter{tocdepth}{2}


\section{Introduction}
The importance of study of the spectrum of self-adjoint operators is well-known. It is equally important to study how the spectrum changes
when a self-adjoint operator, $H_0,$ undergoes a perturbation. Assuming that the essential spectrum, $\sigma_{ess}(H_0),$ stays stable, one important question is to study the net number of eigenvalues
of~$H_0$ which cross a given point $\lambda \notin \sigma_{ess}(H_0)$ as~$H_0$ gets perturbed to $H_1$ via a continuous deformation $H_r, r \in [0,1].$ The resulting number is called the \emph{spectral flow}. 
There are different approaches to the study of spectral flow among which are an intersection number \cite{APS76}, a total Fredholm index \cite{Ph97FIC}, an axiomatic approach \cite{RoSa}
and a recent total resonance index \cite{AzSFnRI}.

Another approach closely related to the approach via the intersection number is as follows. Let~$\clA$ be a real affine space of relatively compact self-adjoint perturbations of~$H_0$ which leave the essential spectrum unchanged, 
and consider the set, $\euR(\lambda),$ of all operators from~$\clA$ for which ~$\lambda$  is an eigenvalue. The set~$\euR(\lambda)$ is an analytic variety, which we call the \emph{resonance set}.
It can be shown that for $\lambda \notin \sigma_{ess}(H_0)$ the variety~$\euR(\lambda)$ has co-dimension 1. As such, this variety divides the affine space~$\clA$ near a point $H_0 \in \clA$ 
into two or more parts, which we will call \emph{resonance cells} or simply \emph{cells}. Assuming for simplicity that~$\lambda$ is a simple eigenvalue of $H_0,$ there will be only two parts, let them be $\euR_+$ and $\euR_-.$ The $\euR_+$ consists of those points in~$\clA$ which have an eigenvalue slightly larger than $\lambda,$ and a similar interpretation applies to~$\euR_-.$ In these terms, if $H_0 \in \euR_-$ and $H_1 \in \euR_+,$ then one should expect the spectral flow of a norm continuous path connecting~$H_0$ to $H_1$ through~$\lambda$ to be $+1.$ A smooth path connecting~$H_0$ and $H_1$ can intersect the resonance set at many points. Thus, the following question arises: let $H_0 \in \clA$ and let $V$ be a self-adjoint operator from the real vector space $\clA_0 := \clA-H_0.$ If~$H_0$ is perturbed in the direction $V,$ where the eigenvalue~$\lambda$ of~$H_0$ will move: to $\euR_+$ 
or $\euR_-?$ What if~$H_0$ is perturbed in the direction $-V?$ Intuitively it is obvious that if $V$ is tangent to the resonance set, then~$H_0$ may get perturbed into the same cell, resulting in zero contribution to the spectral flow through~$\lambda.$ This occurs if the direction $V$ is tangent to the resonance set at $H_0.$

Thus, we conclude that it is interesting to find criteria for the tangency of a direction $V$ at a point~$H_0$ of the resonance set. For a point~$\lambda$ outside the essential spectrum this program was carried out in \cite{AzSFnRI}. The aim of this paper is to do the same for~$\lambda$ inside the essential spectrum. In this case the resonance set~$\euR(\lambda)$ ought to be defined in terms of the limiting absorption principle. We assume that there is a fixed -- rigging -- operator $F$ which is $\abs{H_0}^{1/2}$-compact, and using it we define~$\euR(\lambda)$ as the set of operators $H \in \clA$
for which the norm limit $T_{\lambda + i0}(H)$ does not exist. If $\lambda$ lies outside the essential spectrum then this definition coincides with the one defined above via the eigenvalue equation.
But inside the essential spectrum the character of the set $T_{\lambda + i0}(H)$ changes as it takes into account not only pure point spectrum but also singularly continuous spectrum too. 

This paper is a natural continuation of the study of the essential spectrum initiated in \cite{AzSFIES} and \cite{AzSFnRI}. The reader should consult sections 2 and 3 of \cite{AzSFIES}, which also has a detailed index, for all the relevant definitions which are omitted here. For more motivation for this work one can also consult the upcoming paper \cite{AzDa4}.

As an example, among the results of this paper is the following, see Theorem \ref{T: V-tlV not regular}. Assuming that for a semi-regular point $\lambda$ the limit $T_{\lambda + i0}(H_0)$ does not exist. Are there perturbations $V$ such that none of the limits $T_{\lambda + i0}(H_0+rV),$ $r \in \mbR,$ exist? As it turns out, there are plenty such perturbations $V,$ namely, any regular direction $V$ has a finite rank perturbation with this property.


\section{Tangency properties of directions in the case of~$\lambda \in \sigma_{ess}$}
Suppose, as usual, that a self-adjoint operator~$H_0$ acts on a rigged Hilbert space $(\hilb,F).$ 
One has to assume some sort of compatibility between~$H_0$ and $F,$ which makes things work. To this end, we assume that $F$ is $\abs{H_0}^{1/2}$-compact. 
Let~$\lambda$ be a point inside 
the essential spectrum~$\sigma_{ess}$ of~$H_0$ such that the norm limit $T_{\lambda+i0}(H_0)$ of the sandwiched resolvent 
$$
    T_{\lambda+iy}(H_0) = F R_{\lambda+iy}(H_0)F^*, \quad y>0,
$$
exists. As usual, let~$V$ be a self-adjoint operator from the real Banach space $\clA_0 = F^* \clB_{sa}(\clK)F$ and let $H = H_0+V.$
Elements of $\clA_0$ we call \emph{directions}.
The set of directions~$V$ for which $T_{\lambda+i0}(H_0 + V)$ does not exist we call the \emph{resonance set} and denote it $\euR(\lambda).$ 
The set~$\euR(\lambda)$ is the essential spectrum case analogue of the set of operators for which $\lambda$ is an eigenvalue. 
It is not difficult to show that 
the set~$\euR(\lambda)$ is an analytic variety, in the sense that its intersection with any finite dimensional subspace of $\clA_0$ is the set of zeros of a finite system of real analytic functions. 
For~$\lambda$ outside the essential spectrum, the resonance set has co-dimension one, see \cite{AzSFnRI}. However, inside the essential spectrum 
the co-dimension is usually larger than one. The co-dimension of the resonance set is closely related to the question of path-independence of singular SSF, or in other terms to the question 
of exactness of the infinitesimal singular spectral shift 1-form. 

Since the resonance set is a smooth variety, it makes sense to ask whether a given direction~$V$ is tangent to the resonance set at a given point~$H_0.$ In \cite{AzSFnRI} it was found that for~$\lambda$ outside the essential spectrum a direction~$V$ is tangent if and only if the equation 
$$
   [(1- r R_{\lambda}(H_r)V)]^2 \phi = 0.
$$
has a non-zero solution for some, and therefore, for any non-resonant value of $r$. Here we discuss similar results for~$\lambda $ inside~$\sigma_{ess}.$

\medskip 
One of the results of \cite{AzSFnRI} asserts that for a point~$\lambda$ outside the essential spectrum~$\sigma_{ess}$ of~$H_0$ the order of a regular direction~$V$
is equal to the order of tangency of~$V$ to the resonance set~$\euR(\lambda).$ 
In this paper we show that inside the essential spectrum a tangent to order $k$ direction has order at least~$k,$
but that the reverse in general is not true, in contrast to the case of $\lambda \notin \sigma_{ess}.$ 
The method of proof used in \cite{AzSFnRI} does not work inside the essential spectrum, since it relies on the eigenvalue equation 
$H_s \phi(s) = \lambda \phi(s)$
for an analytic path of operators $H_s,$ which is not available inside~$\sigma_{ess}.$ However, the eigenvalue equation can be rewritten as 
$$(1 + (s-r) R_{\lambda}(H_r)V) \phi(s) = 0.$$ This equation can be adapted for~$\lambda$ inside the essential 
spectrum by writing $$(1 + (s-r) T_{\lambda+i0}(H_r)J) u(s) = 0.$$ It turns out that this little trick allows to overcome difficulties. 

\bigskip
Recall that a real number $\lambda$ is called \emph{semi-regular} or \emph{essentially regular}, if for some operator $H$ from~$\clA$ the limit $T_{\lambda+i0}(H)$ exists. 
Let~$\lambda$ be a semi-regular point. 
Let $H(s)$ be an analytic path in~$\clA(F)$ which consists of resonant at~$\lambda$ operators. Recall that ``$H$ is resonant at~$\lambda$'' means 
$T_{\lambda+i0}(H)$ does not exist, that is, $H \in \euR(\lambda).$ Let~$H_0$ be an operator regular at~$\lambda,$ and $H(s) -H_0 = V(s) = F^*J(s)F.$
Since~$H_0$ is regular at~$\lambda,$ the limit $T_{\lambda+i0}(H_0)$ is defined, and since $H(s)$ is~$\lambda$-resonant, 
there exists a smooth path of vectors $u(s)$ such that 
\begin{equation} \label{F: [1+TJ(s)]u(s)=0}
    \SqBrs{1 + T_{\lambda+i0}(H_0)J(s)}u(s) = 0.
\end{equation}
Now let's take a straight line $H_s = H_0+sV$ which intersects the curve $H(s)$ 
at $H(r_\lambda).$ 

\begin{thm}
If a direction~$V$ is tangent to the~$\lambda$-resonance set then 
$u'(r_\lambda)$ is a resonance vector of order $2$ and, in particular, 
the order of~$V$ is at least~$2.$
\end{thm}
\begin{proof}
By assumption there exists a~$\lambda$-resonant curve $H(s)$ of operators such that~$V$ is tangent to the curve at the point $H(r_\lambda).$
There is a straight line $H_s = H_0+sV$ so that $H_1  = H(r_\lambda),$ $V = H'(r_\lambda),$
and~$H_0$ is regular at $\lambda.$ 
Let $H(s)-H_0=F^*J(s)F.$
Since $H(s)$ is resonant at~$\lambda$ for all $s,$ 
there exists a smooth path of vectors $u(s)$ with $u(r_\lambda)\neq 0$ such that 
(\ref{F: [1+TJ(s)]u(s)=0}) holds. 
We differentiate this equation and take $s = r_\lambda\colon$
$$
    \SqBrs{1 + T_{\lambda+i0}(H_0)J(r_\lambda)}u'(r_\lambda) = - T_{\lambda+i0}(H_0)J'(r_\lambda) u(r_\lambda).
$$
Since $H(r_\lambda) = H_0 + V$ is resonant at~$\lambda$ and $J'(r_\lambda) = V,$ the right hand side of the last display equals~$u(r_\lambda).$ Indeed,
$$
   - T_{\lambda+i0}(H_0)J'(r_\lambda) u(r_\lambda) = - T_{\lambda+i0}(H_0)V u(r_\lambda) = u(r_\lambda),
$$
where the last equality is a resonance equation of order~$1.$ 
Hence, applying the operator in the square brackets to both sides again gives 
$$
    \SqBrs{1 + T_{\lambda+i0}(H_0)J(r_\lambda)}^2 u'(r_\lambda) = 0.
$$
Therefore, $u'(r_\lambda)$ is a resonance vector of order~$2$
and thus the direction $V=F^*J'(r_\lambda)F$ has order at least 2 at $H(r_\lambda).$
\end{proof}

In the proof of the following theorem, for a rigging $F$, we will use the equivalence
\begin{equation} \label{F: FJF=0 then J=0}
     F^*JF = 0 \ \ \text{iff} \ \  J=0.
\end{equation}

\begin{thm}
If a direction~$V$ is tangent to the resonance set to order $k$ then
the vectors $u'(r_\lambda), \ldots, u^{(k-1)}(r_\lambda)$ are resonance vectors of order, respectively, $2, \ldots, k,$ and, in particular,
the order of~$V$ is at least~$k.$
\end{thm}
\begin{proof} 
By premise, there exists a resonant path $H(s),$ such that $H(r_\lambda)$ is the operator at which~$V$ is tangent to the resonance set 
and 
\begin{equation} \label{F: ders of H(s) are 0}
    H'(r_\lambda) = V, \ H''(r_\lambda) = 0, \ \ldots, \ H^{(k-1)}(r_\lambda) = 0.
\end{equation}
Choose a straight line $H_s = H_0+sV$ so that $H_1 = H(r_\lambda),$ and~$\lambda$ is regular at~$H_0.$
We proceed by induction on $k.$  The induction base is the previous theorem. Assume the claim for $k-1.$ Since the path $H(s)$ is resonant, there exists an analytic path $u(s)$ in $\clK$
such that (\ref{F: [1+TJ(s)]u(s)=0}) holds. 
Differentiating this equality $k-1$ times and replacing $s = r_\lambda$ gives, taking into account \eqref{F: ders of H(s) are 0} and~\eqref{F: FJF=0 then J=0},
$$
    \SqBrs{1 + T_{\lambda+i0}(H_0)J(r_\lambda)}u^{(k-1)}(r_\lambda) = - (k-1) T_{\lambda+i0}(H_0) J'(r_\lambda)  u^{(k-2)}(r_\lambda).  
$$
Here we note that $J'(r_\lambda) = J$ where~$J$ is from $V = F^*JF.$
The operator $T_{\lambda+i0}(H_0) J$ preserves the order and the operator $\SqBrs{1 + T_{\lambda+i0}(H_0)J(r_\lambda)}$
decreases it. Thus, since by induction assumption the vector $u^{(k-2)}(r_\lambda)$ has order $k-1,$ 
it follows from the last equality that the vector $u^{(k-1)}(r_\lambda)$ has order $k.$ In particular, the direction~$V$ has order at least~$k.$ 

Further, similar to \cite{AzSFnRI}, one can show that 
$$
   \bfA_{\lambda+i0}(r_\lambda) u^{(k-1)}(r_\lambda) = (k-1)u^{(k-2)}(r_\lambda).
$$
Proof of this equality is this: in the penultimate display replace~$H_0$ by $H_s$ with complex variable $s$ and then take the integral of the both sides
along a small contour which encloses $r_\lambda.$
\end{proof}

Thus, tangent directions have order at least $2.$ However, inside the essential spectrum there are plenty of transversal directions which also have order at least~$2.$
We will call such directions \emph{collar directions}. 

%

Outside the essential spectrum a direction has order $\geq 2$ iff it is tangent to the resonance set.
Therefore, a linear combination of two directions of order $\geq 2$ is also a direction of order $\geq 2.$ 
It is reasonable to believe that there should be an algebraic proof of this fact which hopefully would work 
inside essential spectrum as well.

%

\section{Reduction of direction}
Here we will adapt the results of \cite[\S 7.1]{AzSFnRI} to the case where a spectral point~$\lambda$ belongs to~$\sigma_{ess}.$
The fact that $\lambda\in\sigma_{ess}$ creates some difficulties which we aim to overcome here. 

So, let $\lambda \in \sigma_{ess}$ be a semi-regular point for a self-adjoint $H_0 \in \clA.$
Let $V = F^* J F$ be a regular direction with property~$S,$ see \cite[\S 13.3]{AzSFIES} for the definition of the latter. 
The property~$S$ can be equivalently characterised in many ways \cite[Proposition 13.3.1]{AzSFIES}, one of them is 
$Q_{\lambda-i0} J P_{\lambda+i0} = J P_{\lambda+i0}.$ In particular, a point~$\lambda$ has the property $S$ if and only if the bounded operator $J P_{\lambda+i0}$ 
on $\clK$ is self-adjoint and thus can be treated as a direction. 

Before proceeding further, we note that the property $S$ is a generic property. Suffices to say that all directions of order $1,$ all positive directions
and all directions in case $\lambda \notin \sigma_{ess}$ possess this property. In fact, it is not easy to present examples of directions without the property~$S,$
though this is mainly due to the fact that it is not easy to present examples of points of order~$>1.$

\smallskip 
In this section we will use the following notation.~$H_0$ is semi-regular at~$\lambda,$ $V=F^*JF$ is a regular direction 
and 
$$
    \tlV := F^*\tlJ F, \quad \text{where} \ \ \tlJ:=J P_{\lambda+i0}(H_0,V).
$$
We often use notation $P_+ := P_{\lambda+i0}(H_0,V)$ and $T_+:=T_{\lambda+i0}(H_r),$
where $r$ is some regular point the choice of which should be clear from the context. 
We consider $J \mapsto \tlJ$ as a map on regular elements of $\clA_0.$  
As usual $H_r = H_0+r V,$ and $\tilde H_r = H_0 + r \tilde V.$

The following theorem is the analogue of \cite[Theorem 7.1.1]{AzSFnRI}. It is not much surprising since it is unlikely for a direction at a semi-regular point
not to be regular, but alas in mathematics we should take care of all possibilities even if they are extremely unlikely. 
\begin{thm} \label{T: V reg then tilde V reg} Let~$V$ be a regular direction with property $S$ at a semi-regular point~$\lambda$ for~$H_0.$ Then the direction $\tilde V$ 
is also regular. 
\end{thm}
\begin{proof} By definition, that~$V$ is regular means that $T_{\lambda+i0}(H_r)$ exists for some $r \in \mbR,$ and we have to show the same for 
$T_{\lambda+i0}(\tilde H_r).$
The second resolvent identity applied to $\tilde H_r = H_r + r(\tilde V - V)$ gives 
$$
   T_z(\tilde H_r) =  T_z(H_r + rF^*(\tilde J -J)F) =  \SqBrs{ 1 + r T_z(H_r) (\tilde J-J)  }^{-1}  T_z(H_r).
$$
Thus, $T_{\lambda+i0}(\tilde H_r)$ exists iff the operator $1 + r T_{\lambda+i0}(H_r) (\tilde J-J) $ is invertible. 
Assume the contrary, that is (as the operator is Fredholm with zero index), for some non-zero $\phi$ 
$$
    [1 + r T_{\lambda+i0}(H_r) (\tilde J-J)] \phi = 0.
$$
Since $P_+$ and  $T_{\lambda+i0}(H_r)J$ commute, applying $P_+$ to both sides of the equation above gives
$P_+\phi = 0.$ This gives $\tilde J \phi = 0,$ which combined with the last display gives 
$[1 - r T_{\lambda+i0}(H_r) J] \phi = 0.$ This implies that $\phi \in \im P_+,$ and thus, 
$\phi = P_+ \phi = 0.$
\end{proof}

\begin{thm} \label{T: T(JP+)JP+=T(J)JP+} 
Let $V=F^*JF$ be a regular direction with property $S$ at a semi-regular point~$\lambda$ for~$H_0.$ Then for any non-resonance $s$ 
(w.r.t. $\lambda+i0$) 
$$
   T_{\lambda+i0}(\tilde H_s) JP_+ =    T_{\lambda+i0}(H_s) JP_+.
$$
\end{thm}
\begin{proof} Given Theorem~\ref{T: V reg then tilde V reg}, the proof follows verbatim that of \cite[Theorem 7.1.2]{AzSFnRI} with some notational changes. 
Still, we give it here for reader's convenience. 
We let $A_+(s) = T_{\lambda+i0}(H_s) J.$
By the second resolvent identity, we have 
\begin{equation*}
  \begin{split}
      (E) := T_{\lambda+i0}(\tilde H_s) JP_+  & =  \SqBrs{1 - sT_{\lambda+i0}(H_s)(J-\tlJ) }^{-1} T_{\lambda+i0}(H_s)JP_+ \\
         & = \SqBrs{1 - sA_+(s)(1-P_+) }^{-1} A_+(s)P_+. \\
  \end{split}
\end{equation*}
The operator $\tilde A_+(s) = A_+(s)(1-P_+)$ is the holomorphic part of the Laurent expansion of $A_+(s)$ at $s = 0.$
So, we have for small enough $s$ 
\begin{equation*}
  \begin{split}
      (E)  & = \SqBrs{1 - s \tilde A_+(s)}^{-1} A_+(s)P_+ \\
          & = \SqBrs{1 + s \tilde A_+(s) + s^2 \tilde A^2_+(s) + \ldots }  A_+(s)P_+. 
  \end{split}
\end{equation*}
Since $[ A_+(s), P_+]=0$ and $ \tilde A_+(s)P_+=0$ it follows that $(E) =  A_+(s)P_+,$ as required. 
By analytic continuation, the equality holds for all not necessarily small $s.$
\end{proof}
\begin{thm} \label{T: P(V)=P(tilde V) and }
Let~$V$ be a regular direction with property $S$ at a semi-regular point~$\lambda$ for~$H_0.$ Then 
$$
   P_{\lambda+i0}(H_0,\tilde V) =    P_{\lambda+i0}(H_0,V)
$$
and 
$$
   \bfA_{\lambda+i0}(H_0,\tilde V) =    \bfA_{\lambda+i0}(H_0,V).
$$
\end{thm}
\begin{proof} Given Theorems~\ref{T: V reg then tilde V reg} and~\ref{T: T(JP+)JP+=T(J)JP+}, this proof follows verbatim 
that of \cite[Theorem 7.1.3]{AzSFnRI} with some obvious notational changes. Still, we give this proof as well. 
We prove the second equality, the first one is proved by the same argument. Using the definition of $\bfA_+$ and Theorem~\ref{T: T(JP+)JP+=T(J)JP+},
we have 
\begin{equation*}
  \begin{split}
       \bfA_+(H_0,\tlV) & = \frac 1 {2\pi i} \oint_{C(0)} s T_+(\tilde H_s) \tlJ\,ds \\
       & = \frac 1 {2\pi i} \oint_{C(0)} s T_+(H_s) JP_+\,ds \\
       &  = \bfA_+ P_+  \\
       &  = \bfA_+.
  \end{split}
\end{equation*}

\end{proof}

The following theorem is a direct consequence of the previous theorems (see also Theorems~7.1.4 and~7.1.5 in \cite{AzSFnRI}).
\begin{thm} Let~$V$ be a regular direction with property $S$ at a semi-regular point~$\lambda$ for~$H_0.$ Then 
the resonance matrices of the directions~$V$ and $\tilde V$ are equal, and therefore so are their resonance indices. 
\end{thm} 
Indeed, the resonance matrix of~$V$ is $JP_+(H_0,V)$ and Theorem~\ref{T: P(V)=P(tilde V) and } gives 
$$
   JP_+(H_0,V) = JP_+(H_0,V) \cdot P_+(H_0,V) = \tilde J \cdot P_+(H_0,\tilde V),
$$
where the last operator is the resonance matrix of $\tilde V.$ Further, the resonance index is the signature of the resonance matrix, see \cite[Theorem 9.2.1]{AzSFIES},
and thus the resonance indices of~$V$ and $\tilde V$ are also equal. 

\medskip
Definition of plain homotopic directions is the same as \cite[Definition 7.1.6]{AzSFnRI}.
\begin{thm} Let~$V$ be a regular direction with property $S$ at a semi-regular point~$\lambda$ for~$H_0.$ Then 
the directions~$V$ and $\tilde V$ are plain homotopic.
\end{thm}
\begin{proof} The proof follows verbatim that of \cite[Theorem 7.1.7]{AzSFnRI} with some obvious notational changes
and one more change: everywhere the statement of this form ``$H-\lambda$ is invertible'' should be replaced by 
``$T_{\lambda+i0}(H)$ exists''.
\end{proof}

\begin{lemma} \label{L: P same for all r}
The idempotent $P_z(H_0,rV)$ does not depend on $r\neq 0.$ 
\end{lemma}
\begin{proof}
By definition,
$$
   P_z(H_0,rV) = \frac 1{2\pi i} \oint_{C(0)} T_z(H_{sr})rJ\,ds.
$$
From this one can see that this idempotent does not depend on $r$ including its sign.
Indeed, scaling of $r$ results in scaling of $C(0),$ which does not affect the contour integral (for large $r$ one can always choose $C(0)$
to be small enough). Replacing $r$ by $-r$ also does not change the left side:
\begin{equation*}
  \begin{split}
      P_z(H_0,-V)  & = \frac 1{2\pi i} \oint_{C(0)} T_z(H_0-sV)(-J)\,ds \\
            & = \frac 1{2\pi i} \oint_{-C(0)} T_z(H_0+tV)J\,dt = P_z(H_0,V).
   \end{split}
\end{equation*}
since $t = - s$ also traces the contour $-C(0)$ in the counterclockwise direction.
\end{proof}

At the same time, it is obvious that 
$$
    \ind_{res}(\lambda; H_0,-V) =  -\ind_{res}(\lambda; H_0,V).
$$
Recall that for a semi-regular point~$\lambda$ the resonance set $$\euR(\lambda) := \set{H \in \clA \colon H \text{ is not regular at } \lambda}$$
is an analytic variety of co-dimension $1$ outside the essential spectrum and $>1$ in the inside.
Therefore, inside the essential spectrum we can find an intersection of the resonance set with three dimensional real affine plane
so that the intersection is one-dimensional.
In this 3D section we can transversally deform~$V$ changing the sign of the resonance index. That is, the resonance index
is not homotopically invariant for~$\lambda$ inside~$\sigma_{ess}.$ 
For a semi-regular point of geometric multiplicity $1$ the resonance index can have only the values $\pm 1$ or $0,$ according to the U-turn inequality, see \cite{AzSFIES}.
A transversal direction~$V$ of resonance index $+1$ can be transversally deformed to~$-V$ with resonance index $-1.$ Thus, along the way 
we get a transversal direction of resonance index~$0.$ This indicates that there is a set of ``invisible'' directions of resonance index zero,
which form a barrier for deforming $+1$ resonance index directions to $-1$ resonance index directions. The U-turn inequality shows that these zero resonance index directions
should have the algebraic multiplicity at least~$2.$ Thus, we conclude that 
\begin{prop} Inside the essential spectrum there exist transversal directions of order $\geq 2$ (that is, collar directions).
\end{prop}
This is in contrast to \cite[Theorem 4.3.3]{AzSFnRI} which asserts that outside the essential spectrum
a regular direction has order 1 iff it is transversal. 


\begin{thm} \label{T: T+(sV+t tlV) = T+(sV+tV)} Let $V=F^*JF$ be a regular direction with property $S$ at a semi-regular point~$\lambda$ for~$H_0.$ 
Let $s$ and $t$ be real numbers (such that $s+t\neq 0$) such that $s V + t \tilde V$ is regular.
Then 
$$
   T_{\lambda+i0}(H_0+s V +  t \tilde V) JP_+ =    T_{\lambda+i0}(H_0+(s+t)V) JP_+.
$$
\end{thm}
\begin{proof}
Let $H_{s,t} = H_0+s V +  t \tilde V.$ 
In this notation Theorem~\ref{T: T(JP+)JP+=T(J)JP+} asserts $$T_{\lambda+i0}(H_{s,0})JP_+ = T_{\lambda+i0}(H_{0,s})JP_+.$$
Thus, using alternately the second resolvent identity and the last display twice we get
\begin{equation*}
  \begin{split}
    T_{\lambda+i0}(H_{s,t}) JP_+ & = \SqBrs{1 +t T_{\lambda+i0}(H_{s,0})JP_+ }^{-1}T_{\lambda+i0}(H_{s,0})JP_+ \\
    & = \SqBrs{1 +t T_{\lambda+i0}(H_{0,s})JP_+ }^{-1}T_{\lambda+i0}(H_{0,s}) JP_+\\
    & = T_{\lambda+i0}(H_{0,s+t}) JP_+\\
    & = T_{\lambda+i0}(H_{s+t,0}) JP_+.
  \end{split}
\end{equation*}
\end{proof}

\begin{thm} \label{T: V-tlV not regular} 
For any regular direction~$V$ the direction $V -\tlV$ is not regular. 
\end{thm}
\begin{proof}
That~$H_0$ is semi-regular and non-regular at~$\lambda,$ and~$V$ is a regular direction means that 
the equation 
$$
   (1 - T_+(H_0+V)J)u = 0
$$
makes sense and has a non-zero solution.
To prove the claim that $H_0+V-\tlV$ is semi-regular but not regular at~$\lambda,$
we need to show that the equation 
$$
   (1 - T_+(H_0+2V-\tlV)J)u = 0
$$
also has a non-zero solution. In fact, the same $u$ is a solution of this equation as well. 
Indeed, since $u$ is a solution of the first equation, we have $u = P_+u,$ and therefore, 
by Theorem~\ref{T: T+(sV+t tlV) = T+(sV+tV)}, 
\begin{equation*}
  \begin{split}
   (1 - T_+(H_0+2V-\tlV)J)u  & =    (1 - T_+(H_0+2V-\tlV)J)P_+u  
       \\ & =  (1 - T_+(H_0+V)J)P_+u  \\
       & = 0.
   \end{split}
\end{equation*}   
\end{proof}
The argument of this proof shows that for any $k=1,2,\ldots$ 
$$
    \Upsilon^k_{\lambda+i0}(H_0,V) \subset     \Upsilon^k_{\lambda+i0}(H_0+V-\tilde V,V)
$$

\begin{cor} (a) The resonance set contains a lot of straight lines.\\
(b) At any semi-regular~$H_0$ there are plenty of non-regular directions. Namely, any regular direction has a finite-rank perturbation which is not regular. 
\end{cor}
Outside the essential spectrum ($\lambda \notin \sigma_{ess}$) this corollary can be proved in another simpler way,
and if in addition~$\lambda$ is a simple eigenvalue and~$V$ has order $1$ then a proof becomes particularly trivial.

\bigskip
\begin{thm} For~$H_0$ semi-regular at~$\lambda,$ and~$V$ a regular direction, we have 
\begin{equation*}
  \begin{split}
        P_{\lambda+i0}(H_0,V)  P_{\lambda+i0}(H_0+V-\tilde V,V)  &  =  P_{\lambda+i0}(H_0+V-\tilde V,V) P_{\lambda+i0}(H_0,V) \\
        & = P_{\lambda+i0}(H_0,V)
   \end{split}
\end{equation*}   

and 
\begin{equation*}
  \begin{split}
       P_{\lambda+i0}(H_0,V)  \bfA_{\lambda+i0}(H_0+V-\tilde V,V)  &  =  \bfA_{\lambda+i0}(H_0+V-\tilde V,V) P_{\lambda+i0}(H_0,V) \\
         & = \bfA_{\lambda+i0}(H_0,V)
   \end{split}
\end{equation*}   
\end{thm}
\begin{proof} We have, using Theorem~\ref{T: T+(sV+t tlV) = T+(sV+tV)},
\begin{equation}
   \begin{split}
          2\pi i P_{\lambda+i0}(H_0+V-\tilde V,V) P_+ & = \oint_{C(0)} T_{\lambda+i0}(H_0+V-\tilde V +sV) JP_+\,ds \\
           & = \oint_{C(0)} T_{\lambda+i0}(H_0 +sV) JP_+\,ds = 2\pi i P_+^2 = 2\pi i P_+. \\
   \end{split}
\end{equation}
Thus, 
$$ 
    P_{\lambda+i0}(H_0+V-\tilde V,V) P_+ = P_+.
$$
A display in the proof of Theorem~\ref{T: T+(sV+t tlV) = T+(sV+tV)} implies that $P_+$ commutes with $T_+(H_{s,t})J.$
Therefore, the first display of this proof shows that $P_+$ commutes with     $P_{\lambda+i0}(H_0+V-\tilde V,V).$

Proof of other equalities is similar.
\end{proof}

\begin{thm} \label{T: P+(tilde V + V)=P_+} Under the premise of Theorem~\ref{T: T+(sV+t tlV) = T+(sV+tV)}, for any real numbers $\alpha$ and $\beta,$
such that $\alpha+\beta \neq 0,$ 
$$
   P_{\lambda+i0}(H_0,\alpha V + \beta \tilde V) =    P_{\lambda+i0}(H_0,V)
$$
and,
$$
   \bfA_{\lambda+i0}(H_0,\alpha V + \beta \tilde V) =   \frac 1{\alpha+\beta} \bfA_{\lambda+i0}(H_0,V).
$$
\end{thm}
\begin{proof}
We have, by definition,  
\begin{equation*}
  \begin{split}
       P_{\lambda+i0}(H_0,\alpha V + \beta \tilde V) &  =  \frac 1{2\pi i} \oint_{C}     T_{\lambda+i0}(H_0+s(\alpha V + \beta \tilde V))(\alpha J + \beta \tilde J)\,ds.
  \end{split}
\end{equation*}
We split the last integral into the sum of two integrals and calculate them separately. 
By Theorem~\ref{T: T+(sV+t tlV) = T+(sV+tV)}, for the second summand we have 
$$
    \frac 1{2\pi i} \oint_{C}     T_{\lambda+i0}(H_0+s(\alpha V + \beta \tilde V))\beta \tilde J\,ds = \frac 1{2\pi i} \oint_{C}     T_{\lambda+i0}(H_0+s(\alpha + \beta) V) \beta \tilde J\,ds
      = \frac \beta{\alpha+\beta}P_+.
$$
Now, we consider the first summand
$$
    (E) := \frac 1{2\pi i} \oint_{C}     T_{\lambda+i0}(H_0+s(\alpha V + \beta \tilde V) )\alpha J\,ds.
$$
We have, for small enough $\beta,$ using the fact that $A_+ := A_{\lambda+i0}(s\alpha)$ and $P_+$ commute
\begin{equation*}
  \begin{split}
      T_{\lambda+i0}(H_0+s(\alpha V + \beta \tilde V) )J & =  \SqBrs{1+s\beta  T_{\lambda+i0}(H_0+s\alpha V )JP_+}^{-1}  T_{\lambda+i0}(H_0+s\alpha V )J \\
      & =  \SqBrs{1+s\beta  A_{\lambda+i0}(s\alpha )P_+}^{-1}  A_{\lambda+i0}(s\alpha) \\
      & =  \SqBrs{1 - s\beta  A_+ P_+ + s^2\beta^2  A^2_+P_+ - \ldots}  A_+ \\
      & =  \SqBrs{1 - P_+ + P_+  - s\beta  A_+P_+ + s^2\beta^2  A_+^2P_+ - \ldots}  A_+ \\
      & = \tilde A_+  + \SqBrs{1  - s\beta  A_+ + s^2\beta^2  A_+^2 - \ldots}  A_+ P_+\\
      & = \tilde A_+  + \SqBrs{1  + s\beta  A_+}^{-1}  A_+ P_+\\
      & = \tilde A_+  + A_{\lambda+i0}(s\alpha +s\beta) P_+.\\
  \end{split}
\end{equation*}
The function $\tilde A_+$ is holomorphic at $s = 0,$ so its integral vanishes. Thus,
$$
    (E) = \frac \alpha{\alpha+\beta}P_+.
$$
Combining this with the second display of the proof completes the proof of the first equality, for small enough $\beta.$
For other $\beta$ the equality holds by analytic continuation. 
The second one is proved by the same argument.
\end{proof}
Of course, in these theorems one can replace $\lambda+i0$ by $\lambda-i0.$

%
%
%

\subsection{The map $V \mapsto \tilde V$}
Let~$H_0$ be semi-regular at $\lambda.$ To any regular direction~$V$ with property~$S$ we can assign the direction $\tilde V = F^* \tilde JF = F^*J P_{\lambda+i0}(H_0,V)F.$
We summarise some properties of this map. 

\begin{thm} 
Fix a semi-regular point~$\lambda$ for $H_0 \in \clA.$
The map $V \mapsto \tilde V,$ defined on regular directions $V,$ has the following properties:
\begin{enumerate}
  \item The direction $\tilde V$ is regular.
  \item The direction $V - \tilde V$ is not regular.  
  \item $P_{\lambda+i0}(H_0,\tilde V) = P_{\lambda+i0}(H_0,V).$
  \item $\bfA_{\lambda+i0}(H_0,\tilde V) = \bfA_{\lambda+i0}(H_0,V).$
  \item \label{Item: tilde tilde=tilde} $\tilde {\tilde J} = \tilde J.$  
  \item For any non-zero complex $r$ we have $\tilde{}\, (rJ) = r \tilde J.$
\end{enumerate}
\end{thm}
\begin{proof} All these properties, except (\ref{Item: tilde tilde=tilde}), have been proved before
in Theorems~\ref{T: V reg then tilde V reg}, 
\ref{T: V-tlV not regular},
\ref{T: P(V)=P(tilde V) and } and Lemma~\ref{L: P same for all r}.
 Now we prove (\ref{Item: tilde tilde=tilde}):
$$ 
     \tilde {\tilde J} = \tilde J P_+(H_0,\tilde V) =  JP_+(H_0,V) \cdot P_+(H_0,V) = \tilde J.
$$
\end{proof}

Question: is it true that $\tilde{}\, (J_1+J_2) = \tilde J_1 + \tilde J_2,$ provided $J_1+J_2$ is regular?


\begin{lemma} \label{L: T+ is continuous in V} 
Let~$H_0$ be resonant at $\lambda.$
Let $V \in \clA^\mbC_0$ be a regular direction. 
The operator $T_{\lambda+i0}(H_0+V)$ depends continuously on~$V$ in the $\clA_0$ norm.
\end{lemma}
\begin{proof} Let $V_1 \in \clA_0$ have a small norm. 
Then by the second resolvent identity
$$
   T_{\lambda+i0}(H_0+V+V_1) = T_{\lambda+i0}(H_0+V) \SqBrs{1 + J_1T_{\lambda+i0}(H_0+V)}^{-1}.
$$
Since $J_1$ has small norm, -- by the premise, the operator in the square brackets is invertible and therefore the expression is continuous in $J_1.$
\end{proof}

\begin{thm} \label{T: res index stable for order 1} The resonance index of a direction of order $1$ is stable under small perturbations in $\clA_0$ norm. 
\end{thm}
\begin{proof} Since the resonance index is the signature of the resonance matrix, it is enough to show that the resonance matrix $JP_+(H_0,V)$ depends continuously 
on~$V$ for order $1$ directions and has constant rank. 
The resonance matrix is equal to 
$$
    \frac 1 {2\pi i} \oint_{C(0)} J T_{\lambda+i0}(H_0+sV)J\,ds.
$$
Using compactness of the contour $C(0)$ and Lemma~\ref{L: T+ is continuous in V}, one can show that~$V$ has a neighbourhood in $\clA_0$
such that for all $V'$ from it the operator $T_{\lambda+i0}(H_0+sV')$ exists for all $s \in C(0)$ and is continuous in $s.$
Thus, the integral above is continuous in~$V.$
Since~$V$ has order $1,$ this small perturbation does not generate any new resonance points inside the contour $C(0),$ in addition to the zero resonance point,
and therefore the rank of the resonance matrix is constant. 
Therefore, the integral above with perturbed~$V$ would be the resonance matrix of the perturbed direction, and proof is thus complete. 
\end{proof}



\subsection{The reduction $\tilde V$ outside~$\sigma_{ess}$}
Outside the essential spectrum there is an explicit formula for $P_\lambda,$ see \cite[(5.9.2)]{AzSFnRI},
and therefore there is such a formula for the reduction.
In the case of an order 1 direction~$V$ at a simple point this formula takes very simple form:
$$
    \tilde V = \Scal{\chi,V\chi}^{-1} \Scal{V\chi,\cdot}V\chi,
$$
where $\chi$ is an eigenvector of~$H_0$ corresponding to $\lambda.$ Since the point~$H_0$ is simple, this eigenvector is unique
up to scaling, the choice of which does not affect the formula above. Also, since~$V$ has order $1,$ the inner product $\Scal{\chi,V\chi}$ is non-zero, see \cite{AzSFnRI}.

\subsection{The case of~$V$ with no property $S$}
If a direction~$V$ has no property $S,$ then the operators $Q_-JP_+$ and $Q_+JP_-$ are different, and so we have two candidates 
for a resonance matrix.
This case is still work in progress, but we will present one initial result.

\begin{thm} If~$V$ is a regular direction then each of the directions $Q_-JP_+$ and $Q_+JP_-$ are regular. 
\end{thm}
\begin{proof} Let $\tilde H_r = H_0 + r Q_-JP_+.$ We show that if $T_+(H_r)$ exists then so does $T_+(\tilde H_r).$ We assume, wlog, that $r = 1$
and write $H=H_1$ and $\tilde H = \tilde H_1.$
We have by the second resolvent identity
$$
    T_+(\tilde H) = \SqBrs{1 + T_+(H)(Q_-JP_+-J)}^{-1} T_+(H).
$$
Thus, existence of     $T_+(\tilde H)$ is equivalent to invertibility of $1 + T_+(H)(Q_-JP_+-J).$
Assume it is not invertible. Then for some non-zero $\phi$ we have 
$\SqBrs{1 + T_+(H)(Q_-JP_+-J)} \phi = 0.$

(A) It suffices to show that $P_+ \phi = 0.$ Indeed, if this holds then $$0 = \SqBrs{1 + T_+(H)(Q_-JP_+-J)} \phi = \SqBrs{1 - T_+(H)J} \phi.$$
This equality means that $\phi \in \im P_+,$ and so $\phi = 0.$

(B) We show that $P_+ \phi = 0.$ Apply $P_+$ to the equation $$\SqBrs{1 + T_+(H)(Q_-JP_+-J)} \phi = 0.$$
This gives (using standard properties of $T_+,$ $P_\pm$ and $Q_\pm,$ see \cite[Section 3]{AzSFIES})
$$
    - P_+\phi = (P_+P_-P_+ - P_+) T_+J  P_+ \phi.
$$
The image of $P_+P_-P_+ - P_+$ consists of vectors of type I (see \cite[Lemma 11.2.4]{AzSFIES}). 
Thus $P_+\phi$ is of type I.
But $T_+J$ preserves property of resonance vectors to be of type I. By the same 
\cite[Lemma 11.2.4]{AzSFIES}, the operator $P_+P_-P_+ - P_+$ maps to zero any vector of type I. Hence, from the last display 
 we have $P_+ \phi = 0.$

\medskip For $Q_+JP_-$ proof is the same, --- since $T_+$ exists iff $T_-$ does, all we need to change is to replace $T_+$ by $T_-.$
\end{proof}

\bigskip
{\it Acknowledgements.} The author thanks his wife, Feruza, for financially supporting him during the work on this paper.

\mathsurround 0pt


\end{document}